\documentclass[12pt]{article}
\usepackage{amssymb, latexsym, amsmath, amsthm}
\def\EE{{\mathbb E}}
\def\PP{{\mathbb P}}
\newtheorem{lemma}{Lemma}

\newtheorem{theorem}{Theorem}
\begin{document}
\setlength{\baselineskip}{.25in}
\title{{\bf A Generatingfunctionology Approach to a Problem of Wilf}}
\author{ Pawe\l\  Hitczenko
\\ Department of Mathematics and Computer Science \\ Drexel
University \\ Philadelphia, PA  19104-2875 \\
\texttt{phitczen@mcs.drexel.edu}\\[.1in] Cecil Rousseau
\\ Department of Mathematical Sciences \\ University of Memphis \\
Memphis, TN 38152-3240 \\ \texttt{ccrousse@memphis.edu} \\[.1in]
Carla D. Savage \\ Department of Computer Science \\ North
Carolina State University \\ Raleigh, NC  27695-8206  \\
\texttt{savage@csc.ncsu.edu} }
 \maketitle
\begin{abstract}
  Wilf posed the following problem: determine asymptotically as
  $n \rightarrow \infty$ the probability that a randomly chosen
  part size in a randomly chosen composition of $n$ has multiplicity $m$.
  One solution of this problem has been given by two of the authors \cite{hs}.
  In this paper, we study this question using the techniques of generating
  functions and singularity analysis.
  \end{abstract}
\section{Introduction}
Let $n$ be a positive integer.  A {\em composition} of $n$ with
$p$ parts is a solution of the equation $n = \kappa_1 + \kappa_2 +
\cdots + \kappa_p$ in positive integers $\kappa_1, \kappa_2,
\ldots, \kappa_p$.  We shall write $\kappa = (\kappa_1, \kappa_2,
\ldots, \kappa_p)$ to symbolize the composition.  For example,
there are 16 compositions of 5, namely
\begin{align*}
 & (5)       &  & (4,1)     &  & (1,4)   & & (3,2)  \\
 & (2,3)     &  & (3,1,1)   &  & (1,3,1) & & (1,1,3) \\
 & (2,2,1)   &  & (2,1,2)   &  & (1,2,2) & &  (2,1,1,1) \\
 & (1,2,1,1) &  & (1,1,2,1) &  & (1,1,1,2) & & (1,1,1,1,1).
\end{align*}
The terms $\kappa_1, \ldots, \kappa_p$ are called the {\em parts}
of the composition.   The {\em multiplicity} of a part size is the
number of parts with that size.  For example, in the composition
$(1,1,1,2)$ the multiplicity of 1 is 3 and the multiplicity of 2
is 1.  A {\em partition} of $n$ with $p$ parts is a solution of $n
= \lambda_1 + \lambda_2 + \cdots + \lambda_p$ with $\lambda_1 \geq
\lambda_2 \geq \cdots \geq \lambda_p$.  In \cite{cpsw:rp} Corteel,
Pittel, Savage, and Wilf proved that for every fixed $m \geq 1$,
the probability that a randomly chosen part size in a random
partition of $n$ approaches $1/(m(m+1))$ as $n \rightarrow
\infty$.  Wilf then posed the corresponding problem for
compositions: determine asymptotically (as $n \rightarrow \infty$)
the probability that a randomly chosen part size in a randomly
chosen composition of $n$ has multiplicity $m$.  One solution of
this problem has been given by two of the authors \cite{hs}. In
this note, we address the same question using generating functions
and singularity analysis.

It is well known that there are $2^{n-1}$ compositions of $n$. One
way to arrive at this result uses generating functions.  The
generating function for compositions with $p$ parts is
\[
(z + z^2 + z^3 + \cdots \; )^p = \left( \frac{z}{1-z} \right)^p,
\]
and summing over $p$ we have the generating function for all
compositions:
\[
G(z) = \sum_{p=1}^{\infty} \left( \frac{z}{1-z} \right)^p =
\frac{z}{1-2z}.
\]
The coefficient of $z^n$ in the expansion of $G(z)$, denoted by
$[z^n]G(z)$, is the number of compositions of $n$, and clearly
$[z^n]G(z) = 2^{n-1}$. From an analytic point of view, the fact
that the number of compositions of $n$ is asymptotically (as well
as exactly) $2^{n-1}$ is associated with the fact that the
generating function is a rational function for which the pole
nearest the origin (the only pole in this case) is simple and
located at $z = \frac{1}{2}$.  Our solution of Wilf's problem uses
the same principle.

We shall use the following notation.  The probability of the event
$A$ is denoted by $\PP(A)$, and the expected value of a random
variable $X$ is denoted by $\EE(X)$.  The natural logarithm and
base 2 logarithm are denoted by $\log n$ and $\log_2 n$,
respectively.

 To
state the problem more precisely, suppose that a composition
$\kappa$ is selected uniformly at random from the set of all
$2^{n-1}$ compositions of $n$.  Then out of the set of part sizes
in $\kappa$, a part size $k$ is chosen uniformly at random. Let
$A_n^{(m)}$ denote the event in which $k$ has multiplicity $m$.
For example, inspection of the 16 partitions of 5 shown above
yields \[ \PP(A_5^{(1)}) = \frac{5}{8}, \qquad \PP(A_5^{(2)}) =
\frac{3}{16}, \qquad \PP(A_5^{(3)}) = \frac{1}{8}, \qquad
\PP(A_5^{(5)}) = \frac{1}{16}, \] and otherwise $\PP(A_5^{(m)}) =
0$. The object is to determine $\PP(A_n^{(m)})$ asymptotically as
$n \rightarrow \infty$. We shall find that $\PP(A_n^{(m)})$ tends
to 0 at the rate $1/\log n$.  It then turns out that $\log n \cdot
\PP(A_n^{(m)})$ does not have a limit, but oscillates about the
value $1/m$ as $n \rightarrow \infty$.

\section{Results}
The answer to Wilf's question is given in the following theorem,
first proved in \cite{hs}.
\begin{theorem}
Let $A_n^{(m)}$ be the event in which a randomly selected part
size in a randomly selected composition of $n$ has multiplicity
$m$. Then
\[
\log n \cdot  \PP(A_n^{(m)}) = (1 + o(1))\left( \frac{1}{m} +
F(\{\log_2 n\}) \right), \qquad n \rightarrow \infty,
\]
where $\{a\} = a - \lfloor a \rfloor$ is the fractional part of
$a$ and
\[
F(x) = \frac{2}{m!} \text{Re} \; \sum_{p=1}^{\infty} e^{-2 \pi i p
x} \, \Gamma \left( 1 + i \, \frac{2 \pi p}{\log 2} \right),
\]
with $\Gamma$ denoting the gamma function.

\end{theorem}
Using well-known facts about the gamma function ($\Gamma(1+z) = z
\Gamma(z)$ and $\Gamma(z)\Gamma(1-z) = \pi/\sin(\pi z)$), we
obtain
\[ F(x) = \frac{2}{m!} \sum_{p=1}^{\infty} \left( \frac{p
\alpha}{\sinh(p\alpha)} \right)^{1/2} \cos(2 \pi p x - \phi_p),
\]
where $\alpha = 2 \pi^2/\log 2$ and $\phi_p$ is the argument of
$\Gamma(1 + i \, 2 \pi p/\log 2)$.  This series converges quite
rapidly, and its sum may be approximated by the first term. But
even the first term is quite small since $2(\alpha/\sinh
\alpha)^{1/2} \approx 10^{-5}$.  Thus for large $n$ one finds that
$\log n \cdot \PP(A_n^{(m)})$ is quite close to $1/m$, but there
is a residual dependence on $\{\log_2 n\}$.  In the treatment
given here using generating functions and singularity analysis,
the proof of Theorem 1 will reduce to a well-known calculation
after we have established the appropriate sequence of lemmas.

 Let $\kappa$ be a composition of $n$. Then ${\cal D}(\kappa)$
will denote the set of distinct part sizes in $\kappa$, and ${\cal
M}_m(\kappa)$ will denote the set of part sizes of $\kappa$ that
have multiplicity $m$.
\begin{lemma} 
For a random composition of $n$,
\[
\PP(k \in {\cal M}_m(\kappa)) = \frac{1}{2^{n-1}} [z^n]
\frac{z^{km} (1-z)^{m+1}}{(1 - 2z + z^k(1-z))^{m+1}}.
\]
\end{lemma}
\begin{proof}
Let $G_k(z,w)$ be the two-variable generating function in which
$[z^n w^m] G_k(z,w)$ is the number of compositions of $n$ in which
$k$ has multiplicity $m$. To construct such a generating function,
we first note that the contribution made by compositions with $p$
(not necessarily distinct) parts is
\[
(z + z^2 + \cdots + wz^k + z^{k+1} + \cdots)^p = \left(
\frac{z}{1-z} + (w-1)z^k \right)^p.
\]
Thus
\begin{align*}
G_k(z,w) & =  \sum_{p=1}^{\infty} \left( \frac{z}{1-z} + (w-1)z^k
\right)^p \\[.1in] & = \frac{z + (w-1)z^k(1-z)}{1 - 2z -
(w-1)z^k(1-z)}
\\[.1in] & = \frac{1-z}{1-2z - (w-1)z^k(1-z)} - 1
\\[.1in]
& = \frac{1-z}{1-2z + z^k(1-z) - wz^k(1-z)} - 1.
\end{align*}
Since there are $2^{n-1}$ compositions of $n$, we then have
\begin{align*}
\PP(k \in {\cal M}_m(\kappa)) & = \frac{1}{2^{n-1}} [z^n
w^m]\frac{1-z}{1-2z + z^k(1-z) - wz^k(1-z)} \\[.1in] & =
\frac{1}{2^{n-1}} [z^n] \frac{z^{km}(1-z)^{m+1}}{(1 - 2z +
z^k(1-z))^{m+1}},
\end{align*}
as claimed. \end{proof}
\begin{lemma} 
The polynomial $1 - 2z + z^k(1-z)$ has precisely one zero $z =
\rho$ satisfying $|z| \leq 1$.  This zero is given by
\[
\rho = \frac{1}{2} + \frac{1}{2^{k+2}} + O \left( \frac{k}{2^{2k}}
\right), \qquad k \rightarrow \infty.
\]
For all $k \geq 1$,
\[
\exp \left( - \frac{n}{2^k} \right) < \frac{1}{(2 \rho)^n} < \exp
\left( - \frac{n}{2^{k+2}} \right).
\]
\end{lemma}
\begin{proof} For the first part, simply observe that if
$z = e^{i \theta}$ then $|1-2z|^2 = 5 - 4 \cos \theta$ and
$|1-z|^2 = 2 - 2 \cos \theta$, so $|1-2z| > |1-z|$ for all $z$
with $|z| = 1$.  Apply Rouch\'e's theorem.  To get the approximate
location of $\rho$, write $\rho = \frac{1}{2} + \epsilon$ and
substitute into $1 - 2 \rho + \rho^k(1-\rho) = 0$. This yields
\[
\epsilon = \frac{1}{2^{k+2}} + \frac{k-1}{2^{k+1}} \, \epsilon +
\frac{k(k-3)}{2^k} \, \epsilon^2 + \cdots \; ,
\]
and thus the stated result by iteration.  Next we prove
\[
\frac{1}{2 - 2^{-(k+1)}} < \rho < \frac{1}{2} + \frac{1}{2^{k+1}}.
\]  A simple calculation shows that
$Q(x) = 1 - 2x + x^k(1-x)$ decreases on $(0,1)$. Set $a = 1/(2 -
2^{-(k+1)})$ and $b = \frac{1}{2} + 2^{-(k+1)}$. Then we find that
$Q(a)
> 0$ and $Q(b) < 0$, so $a < \rho < b$.  Then since $(1+x)^n <
\exp(nx)$ for $x > -1$, we have
\[
\frac{1}{(2 \rho)^n} < \left(1 - \frac{1}{2^{k+2}} \right)^n <
\exp \left( - \frac{n}{2^{k+2}} \right) \qquad \text{and} \qquad
\frac{1}{(2\rho)^n} > \frac{1}{(1 + 2^{-k})^n} > \exp \left(-
\frac{n}{2^k} \right).
\]
\end{proof}
We shall show that the number of distinct part sizes $|{\cal
D}(\kappa)|$ of a random composition of $n$ satisfies $|{\cal D}|
\sim \log_2 n$ with probability $1-o(1)$.   The underlying
probabilistic considerations are given in the following lemma.
\begin{lemma} 
Let $X = \sum I_j$ where $(I_j)$ are indicator random variables.
Suppose that $\PP(I_k) = p_{k,n} = p_k$.   If $a$ and $b$ are
chosen so that both $\sum_{j \leq a}(1-p_j)$ and $\sum_{j > b}
p_j$ are $o(1)$, then
\[
\PP(a \leq X \leq b) \geq 1 - o(1).
\]
\end{lemma}
\begin{proof}
For all $a \leq b$ we have
\[
\PP(a \leq X \leq b) = 1 - \PP(\{X < a\} \cup \{X > b\}) \geq 1 -
\PP(X < a) - \PP(X > b).
\]
Now, denoting for simplicity a set and its indicator by the same
symbol,
\[
\PP(X > b) \leq \PP \left(\bigcup_{j > b}I_j \right) \leq \sum_{j
> b} \PP(I_j) = \sum_{j > b} p_j,
\]
and
\[
\PP(X < a) \leq \PP\left(\bigcup_{j \leq a} I_j^c\right) \leq
\sum_{j \leq a}(1 - \PP(I_j)) = \sum_{j \leq a}(1-p_j).
\]
Hence, if $a$ and $b$ are chosen so that both $\sum_{j \leq
a}(1-p_j)$ and $\sum_{j > b} p_j$ are $o(1)$ we get
\[
\PP(a \leq X \leq b) \geq 1 - o(1).
\]
\end{proof}
\begin{lemma} 
Let $\kappa$ be a random composition of $n$. As $n \rightarrow
\infty$ the number of distinct part sizes $|{\cal D}(\kappa)|$
satisfies $|{\cal D}(\kappa)| \sim \log_2 n$ with probability $1 -
o(1)$.
\end{lemma}
\begin{proof} As a special case of Lemma 1, the probability that
$k$ has multiplicity 0 in the random composition $\kappa$ is
\[
\PP(k \in {\cal M}_0(\kappa)) = \frac{1}{2^{n-1}} [z^n]
\frac{1-z}{1-2z + z^k(1-z)}, \qquad n \geq 1.
\]
Hence
\[
\PP(k \in {\cal D}(\kappa)) = 1 - \frac{1}{2^{n-1}} [z^n]
\frac{1-z}{1-2z + z^k(1-z)}.
\]
From Lemma 2, the rational function $(1-z)/(1-2z + z^k(1-z))$ is
analytic for $|z| \leq 1$ except for a simple pole at $z = \rho
\approx \frac{1}{2}$. The residue is $-(1-\rho)/(2 + (k+1)\rho^k -
k \rho^k)$. By standard arguments \cite[\S 5.2]{hw:gf},
\[
[z^n] \frac{1-z}{1-2z + z^k(1-z)} = \left( \frac{1 - \rho}{2 +
(k+1)\rho^k - k \rho^{k-1}} \right) \frac{1}{\rho^{n+1}} + O(1).
\]
By Lemma 2,
\[
 \frac{1 - \rho}{2 + (k+1)\rho^k - k \rho^{k-1}} = \frac{1}{4}\left( 1 + O
 \left( \frac{k}{2^k} \right) \right).
\]
 Hence we have
\[
\frac{1}{2^{n-1}} [z^n] \frac{1-z}{1-2z + z^k(1-z)}  = \frac{1}{(2
\rho)^{n+1}}\left( 1 + O \left( \frac{k}{2^k} \right) \right).
\]
 Using the general bound from Lemma 2
\[
\exp \left( - \frac{n}{2^k} \right) < \frac{1}{(2 \rho)^n} < \exp
\left( - \frac{n}{2^{k+2}} \right),
\]
we see that
\begin{align*}
\PP(k\in {\cal D}(\kappa)) & = 1- \frac{1}{(2\rho)^{n+1}} \left(1+
O\left(\frac{k}{2^{k}} \right)\right)
\\[.1in] & \le 1-\exp\left\{- \frac{n+1}{2^k}\right\}
\left(1+ O\left(\frac{k}{2^{k}} \right)\right)
\\[.1in] & \le
  \frac{n+1}{2^k} + \exp\left\{- \frac{n+1}{2^k}\right\}
\cdot O\left(\frac{k}{2^{k}}\right) ,
 \end{align*}
so that letting $b=\lfloor \log_2n\rfloor +\log\log n$ we get
\[\sum_{k>b}\PP(k\in {\cal D}(\kappa)) =
 O \left(
\frac1{\log n}\right). \] Similarly,
\[
\PP(k \in {\cal D}(\kappa)) \ge
 1-\exp\left\{- \frac{n+1}{2^{k+1}}\right\}
\left(1+ O\left(\frac{k}{2^{k}} \right)\right), \] that is
\[
1-\PP(k\in {\cal D}(\kappa)) \le
 \exp\left\{- \frac{n+1}{2^{k+1}}\right\}
\left(1+ O\left(\frac{k}{2^{k}} \right)\right). \] Consequently,
for any positive $a$,
\begin{align*}
\sum_{1 \le k\le a}\left(1-\PP(k\in{\cal D}(\kappa))\right) & \le
C\sum_{1\le k\le a}\exp\left\{-\frac{n+1}{2^{k+1}}\right\}
 \\[.1in] & =
C\sum_{0\le r < a}\exp\left\{-2^r\frac{n+1}{2^{a+1}}\right\}
\\[.1in] & \le
C \; \sum_{r\ge0}\exp\left\{-(r+1)\frac{n+1}{2^{a+1}}\right\}
\\[.1in] & =
C \, \frac{\exp\left\{-(n+1)/2^{a+2}\right\}}
{1-\exp\left\{-(n+1)/2^{a+2}\right\}},
\end{align*}
and thus \[ \sum_{1 \le k\le a}\left(1-\PP(k\in{\cal
D}(\kappa))\right) = O\left(\frac1{\log n}\right),
\]
provided $a\le\lfloor\log_2n\rfloor-\log\log n$. Hence, by Lemma 3
applied to $I_k=\{k\in{\cal D}(\kappa)\}$, $|D(\kappa)| \sim
\log_2 n$ with probability $1 - o(1)$. \end{proof}

Given a random composition $\kappa$, the probability that a
randomly selected part thereof has multiplicity $m$ is $|{\cal
M}_m(\kappa)|/|{\cal D}(\kappa)|$.  Lemma 4 greatly simplifies the
basic problem. Since so doing amounts to the neglect of a set of
compositions with total probability measure $o(1)$, we may assume
that as $n \rightarrow \infty$ the randomly selected composition
$\kappa$ satisfies $|{\cal D}(\kappa)| \sim \log_2 n$.
  Thus $\PP(A_n^{(m)}) \sim \EE(|{\cal M}_m|)/\log_2 n$.

Now we wish to study the asymptotic behavior of $\PP(k \in {\cal
M}_m(\kappa))$, with the aim of estimating \[ \EE(|{\cal M}_m|) =
\sum_k \PP(k \in {\cal M}_m(\kappa)). \]
\begin{lemma} 
The expected value of $|{\cal M}_m|$ is given by
\[
\EE(|{\cal M}_m|) = (1 + o(1)) \, \frac{n^m}{m!} \sum_k 2^{-km}
\exp(-n/2^k).
\]
\end{lemma}
\begin{proof}
As we found in Lemma 1, the relevant generating function is
\[
G(z) = \frac{1}{2^{n-1}} \frac{P(z)}{Q^{m+1}(z)}, \quad
\text{where} \quad P(z) = z^{km}(1-z)^{m+1}, \; Q(z) =
1-2z+z^k(1-z).
\]
Recall that $Q$ has a simple zero at $z = \rho \approx
\frac{1}{2}$ and no other zeros in $\{z: \; |z| \leq 1\}$.  In a
deleted neighborhood of $\rho$, we have the Laurent expansion
\[
\frac{P(z)}{Q^{m+1}(z)} = \sum_{r = 1}^{m+1}
\frac{c_{-r}}{(z-\rho)^r} + \sum_{s=0}^{\infty} c_s (z - \rho)^s.
\]
The asymptotic behavior of $[z^n] P(z)/Q^{m+1}(z)$ is governed by
the principal part, more specifically by the $r = m+1$ term. In
view of
\[
[z^n] (1-z)^{-\alpha} = \binom{n+\alpha-1}{n},
\]
a simple calculation gives \[ \PP(k \in {\cal M}_m(\kappa)) =
\binom{n+m}{m} \frac{2 P(\rho)}{(- \rho Q'(\rho))^{m+1}} \,
\frac{1}{(2 \rho)^n} \left( 1 + O \left( \frac{1}{n} \right)
\right).
\]
Set
\[
q(n) = \frac{\log n - \log \log n - \log(4(m+1))}{\log 2},
\]
and note that if $k < q(n)$ then $2^k < n/(4(m+1)\log n)$, so
\[ \frac{n^m}{(2 \rho)^n} < n^m \exp \left( -
\frac{n}{4 \cdot 2^k} \right) < n^m \exp(-(m+1)\log n) =
\frac{1}{n}.
\]
Hence we have \[ \sum_{k \leq q(n)} \PP(k \in {\cal M}_m(\kappa))
= O \left( \frac{\log n}{n} \right), \qquad n \rightarrow \infty.
\]

In view of the fact just noted, in estimating $\sum_k \PP(k \in
{\cal M}_m)$, we can limit ourselves to cases where $k > q(n)$. In
that case
\[
\frac{1}{(2 \rho)^n} = \exp \left(- \frac{n}{2^{k+1}} \right)
\left( 1 + O \left( \frac{(\log n)^3}{n^2} \right) \right), \qquad
n \rightarrow \infty.
\]
Now
\[
P(\rho) = 2^{-km}2^{-(m+1)}\left(1 + O \left( \frac{\log n}{n}
\right) \right) \qquad \text{and} \qquad \rho \, Q'(\rho) = -1 + O
\left( \frac{\log n}{n} \right),
\]
so \[ \PP(k \in {\cal M}_m(\kappa)) = \left(1 + O  \left(
\frac{\log n}{n} \right) \right) \, \frac{n^m}{m!}  \, 2^{-(k+1)m}
\exp \left( - \frac{n}{2^{k+1}} \right).
\]
It is now evident that the contribution to the sum $\sum_{k} \PP(k
\in {\cal M}_m(\kappa))$ from those terms with $k > \log_2 n +
\log \log n$ is $o(1)$, so there are  $O(\log \log n)$ terms in
the sum that make a nontrivial contribution.  Thus the bound on
the error for an individual term suffices to give the correct
asymptotic result for the sum.  Replacing $k+1$ by $k$ in the sum,
we have the stated result.
\end{proof}
\begin{proof}[Proof of Theorem 1]
The computational problem that remains is the asymptotic
evaluation of
\[
\frac{n^m}{m!} \sum_{k=1}^{\infty} 2^{-km} \exp(-n/2^k).
\]
Problems of this kind occur frequently in probability theory and
the analysis of algorithms, and now there are different methods
available for their study, and these methods are described in
several excellent references \cite[chapter 7]{fs:ac}. We sketch an
approach due to N. G. de Bruijn, which is described in \cite[pp.
131--134]{kn:ss} and elsewhere.  A special case ($m=1$) of the
above sum is treated in \cite{ccr:322}. The starting point is
Mellin's formula
\[
\exp(-w) = \frac{1}{2 \pi i} \int_{\sigma - i \infty}^{\sigma + i
\infty} w^{-z} \Gamma(z) \, dz, \qquad w, \sigma > 0.
\]
Substituting this representation (with $\sigma = m - \frac{1}{2}$)
for $\exp(-n/2^k)$ and using uniform convergence, one obtains
\[
\sum_{k=1}^{\infty} 2^{-km} \exp(-n/2^k) = \frac{1}{2 \pi i}
\int_{\sigma - i \infty}^{\sigma + i \infty} \frac{n^{-z} \,
\Gamma(z)}{2^{m-z} - 1} \, dz.
\]
Then by the residue theorem, \begin{align*} \frac{n^m}{m!} &
\sum_{k=1}^{\infty} 2^{-km} \exp(-n/2^k) \\[.1in] & = \frac{1}{m!
\log 2} \left\{ (m-1)! + 2 \text{Re} \; \sum_{p=1}^{\infty} e^{-2
\pi i p \log_2 n} \, \Gamma \left(m + i \, \frac{2 \pi p}{\log 2}
\right) \right\}(1 + o(1)).
\end{align*}
The stated result for $\log n \cdot \PP(A_n^{(m)})$ follows.
\end{proof}


\begin{thebibliography}{9}
\bibitem{cpsw:rp} S. Corteel, B. Pittel, C. D. Savage, H. S.
Wilf, On the multiplicity of parts in a random partition, {\em
Random Structures Algorithms} {\bf 14} (1999), pp. 185--197.
\bibitem{fs:ac}
P. Flajolet, R. Sedgewick, {\em Analytic Combinatorics}, to be
published.  In part available at
\texttt{<http://pauillac.inria.fr/algo/flajolet/Publications/publist.html>}.
\bibitem{hs}
P. Hitczenko,  C. D. Savage, On the multiplicity of parts in a
random composition of a large integer, preprint. Available at
\texttt{<http://www.csc.ncsu.edu/faculty/savage/>}.
\bibitem{kn:ss} D. E. Knuth, {\em The Art of Computer Programming. Vol. III:
Sorting and Searching},
Addison-Wesley, Reading, Massachusetts, 1973.
\bibitem{ccr:322} C. C. Rousseau, Solution of Problem 322,
{\em Canad. Math. Bull} {\bf 26} (1983), pp. 375--377.
\bibitem{hw:gf}
H. S. Wilf, {\em Generatingfunctionology}, 2nd ed., Academic
Press, San Diego, 1994. Available at
\texttt{<http://www.cis.upenn.edu/$\sim$wilf/>}.\end{thebibliography}
\end{document}